# Efficient and Precise Calculation of the Confluent Hypergeometric Function

Author: Alan Herschtal[1]

**Abstract**

Kummer's function, also known as the confluent hypergeometric function (CHF), is an important mathematical function, in particular due to its many special cases, which include the Bessel function, the incomplete Gamma function and the error function (erf). The CHF has no closed form expression, but instead is most commonly expressed as an infinite sum of ratios of rising factorials, which makes its precise and efficient calculation challenging. It is a function of three parameters, the first two being the rising factorial base of the numerator and denominator, and the third being a scale parameter. Accurate and efficient calculation for large values of the scale parameter is particularly challenging due to numeric underflow and overflow which easily occur when summing the underlying component terms. This work presents an elegant and precise mathematical algorithm for the calculation of the CHF, which is of particular advantage for large values of the scale parameter. This method massively reduces the number and range of component terms which need to be summed to achieve any required precision, thus obviating the need for the computationally intensive transformations needed by current algorithms.

**Acknowlededements**: The encouragement and guidance of Stephane Heritier, Monash University, School of Public Health and Preventive Medicine, is acknowledged. Proof reading by Liz Ryan, Monash University, School of Public Health and Preventive Medicine, is acknowledged.

[1] Monash University

School of Public Health and Preventive Medicine

Melbourne 3004, Australia

Alan.Herschtal@monash.edu
ORCID: 0009-0004-2061-4787

# 1. Introduction

We commence by describing the Confluent Hypergeometric Function (CHF) and some of its applications. It is a function of 3 parameters, and is commonly denoted as $_1F_1(a;b;z)$ where:

$$_1F_1(a;b;z) = M(a,b,z) = \sum_{n=0}^{\infty} \frac{a^{(n)} z^n}{b^{(n)} n!}$$

$a^{(n)} = a(a+1)(a+2)\cdots(a+n-1)$ is the rising factorial function, a generalisation of the standard factorial where the product starts at $a$ instead of 1. $a^{(0)} = 1$. A discussion of the CHF can be found in Slater [1]. $z$ is commonly referred to as the scale parameter.

The CHF has several important special cases, which have broad application in physics, biology, signal processing and statistics [2]. When $b = 2a$ we have the Bessel function, important in describing wave propagation in a cylinder, in signal processing [3], and in describing the helical structure of DNA. The error function $\text{erf}(x)$ is a scaling of the CHF with $a = \frac{1}{2}$, $b = \frac{3}{2}$, $z = -x^2$. The CHF is closely related to Laguerre polynomials, with application in quantum mechanics, and Hermite polynomials [4] with applications in signal processing, probability and many areas of physics. With $a = 1$, $b = s+1$ and $z = x$, it is also related to the incomplete Gamma function $\gamma(s,x)$. The Poison-Beta distribution, increasingly used in RNA-seq data analysis [5], can be expressed as a scaled reparameterization of the CHF with $b > a$ [6].

Numerical computation of the CHF, like that of many functions without closed form expression, involves a patchwork of solutions, each of which covers a portion of the parameter space [7, 8]. Power series expansions and Buchholz polynomials [9] are efficient and accurate for small $|z|$, but both break down for large $|z|$, ($\gtrsim 100$) [10, 11], with their limit being dependent on the number of bits of floating point precision available. A modification of a power series expansion which requires only a single division can be found in Nardin et at least [12], but this is often less accurate than the power series expansion upon which it is based [13]. Multi-point Taylor expansions increase accuracy in the small to medium $|z|$ part of the parameter space [11], but do not help for large $|z|$. Asymptotic series work well for large $|z|$, but lose accuracy or fail to converge for large $|a|$ or $|b|$. Gauss-Jacobi quadrature [14] based on expression of the CHF as an integral [15] can be effective but is only valid for $b > a$. Additionally, it is computationally very intensive compared with series expansion methods and is thus only recommended when these methods fail. Recurrence relations [16, 17] can be used to transform the large $|a|$ or $|b|$ case into one with $|a|$ or $|b|$ small enough to be managed by one of the aforementioned methods, but incur an additional layer of computational complexity due to the application of the recurrence relation as a pre-processing step. The use of software arithmetic with variable precision and unlimited exponents can be used to circumvent loss of accuracy for large |z|. However this adds a computational overhead of approximately 100-fold compared to machine floating-point arithmetic at low precision[8]. Many of the available methods are summarized by Muller [10] and by Pearson et at least [7, 13].

This work introduces an algorithmically simple, fast and accurate method of calculating the CHF for a very wide range of $a$, $b$, and $z$ values. It is of particular value for large |z|, the most challenging region of the parameter space.

## 2. Method

We note that $M(a, b, z)$ is a polynomial in $z$, and that by using the definition of the rising factorial components $a^{(n)}$ and $b^{(n)}$, any term can be expressed in terms of the previous by a first order recurrence relation. Denoting the $n^{th}$ term in the series by $m_n$ so that $M(a, b, z) = \sum_{n=0}^{\infty} m_n$, we have

$$m_{n+1} = m_n \times \frac{a+n}{b+n} \times \frac{z}{n+1} \qquad (1)$$

Starting with $m_0 = 1$ and applying the recurrence relation above until convergence provides a direct and simple method of computing $M(a, b, z)$. Algorithms to implement this power series approach can be found in Muller [10] and in Relph [18]. It is common to stop summing the component terms $m_n$ when the ratio of the most recent term to the sum of all previous terms falls below a pre-specified tolerance [7]. This recurrence relation approach is the simplest of a variety of series-based approaches developed for calculation of the CHF [11]. However, they all break down for large |z| due to the fact that calculation starts at $n = 0$, and for large |z|, $m_n$ for small $n$ is many orders of magnitude smaller than $m_n$ at its maximum ($\max_n m_n$). Since calculation of $m_n$ at its maximum is dependent on calculation of $m_n$ at smaller values of $n$ via the recurrence relation in Eq. 1, the loss of bits of precision in calculating $M(a, b, z)$ is roughly equal to $\log_2(m_n)$. Even for somewhat moderate values of |z| this loss of precision can override the machine precision such that very few or even no bits of precision remain in the final sum calculated for $M(a, b, z)$.

Further, even when the range of values of $m_n$ to be calculated by the recurrence relationship and then summed remains within machine precision limits, the number of terms that need to be included in the summation despite making a trivial contribution to the overall sum becomes large with increasing |z|, thus reducing computational efficiency.

This work presents an algorithmically straightforward method of pre-computing the range of values of $n$ for which $m_n$ is not trivially small relative to its maximum value, $\max_n m_n$, and hence contributes non-trivially to $M(a,b,z)$. We refer to this range of values as the 'Region of Interest' (RoI). $M(a,b,z)$ may then be calculated by means of a series-based approach, which, rather than commencing at $n = 0$, commences at the lower edge of the RoI. This overcomes both the precision and computational efficiency limitations for large $|z|$ described above. The recurrence relation expressed in Eq. 1, which is the algorithmically simplest method of calculation of the CHF, is ordinarily restricted to small to moderate values of $|z|$. The proposed method extends its use to an almost unlimited range of $|z|$. We focus on the parameter sub-space in which $a, b, z$ are all real, and note that the case where $z < 0$ can easily be transformed to a case where $z > 0$ using Kummer's transformation [19, 20], $M(a,b,z) = e^z M(b-a, b, -z)$, and thus need not be dealt with directly.

Formally, we define the RoI as $\{n: m_n / \max_n m_n > \varepsilon\}$, being the range of values of $n$ outside of which the values of $m_n$ are trivially small relative to $\max_n m_n$, and hence can be ignored in the computation of $M(a,b,z)$. A precision parameter, $\varepsilon$, defining the threshold for what is considered to be 'trivially small', is decided upon a priori and will depend on the application of interest. Typically, one would require at least $\varepsilon \lesssim 10^{-6}$, with much smaller values potentially being required for some applications.

## 2.1. Finding the Maximum CHF Component Term

In order to determine the RoI, we first set about finding the value of $n$ at which $m_n$ is maximised, $n_\mathrm{m} = \mathrm{argmax}_n\, m_n$. We commence by establishing that $m_n$ is unimodal and thus that the region around its mode is the only region of interest that needs to be considered, i.e. that the RoI is a single contiguous region. We denote the ratio of two consecutive terms in the series $m_n$ by $\omega_n$,

$$\omega_n = \frac{m_{n+1}}{m_n} = \frac{a+n}{b+n} \times \frac{z}{n+1} \qquad (2)$$

A necessary and sufficient condition for $m_n$ to be a local maximum is that $\omega_n < 1$ and $\omega_{n-1} > 1$. Although for the purpose of calculating the CHF, $n$ is integer valued, conceptually we may consider it as continuously valued for the purpose of examining the properties of $\omega_n$. Since $\omega_{n_\mathrm{m}-1}$ and $\omega_{n_\mathrm{m}}$ straddle 1, given that $\omega_n$ is defined for all $n$ (except for a discontinuity at $n = -1$ which is not of interest since $n \geq 0$, and another at $n = -b$, which is not of interest unless $b < 0$ for the same reason), there must exist a value of $n$ in $(n_\mathrm{m} - 1, n_\mathrm{m})$ for which $\omega_n = 1$. We thus solve for $n$ in

$$\frac{a+n}{b+n} \times \frac{z}{n+1} = 1$$

This may be expressed as a quadratic $n^2 + (b + 1 - z) \times n + b - az = 0$, the roots of which are

$$n = \frac{(z - 1 - b) \pm \sqrt{(z - 1 - b)^2 + 4(az - b)}}{2} \qquad (3)$$

Only solutions for Eq. 3 for which $n > 0$ are of interest. If no such solution exists then the RoI commences at $n = 0$ and thus there is no computational advantage in identifying the RoI and restricting computation to the region over traditional methods which commence at $n = 0$.

It is straightforward to show that for $z > 1 + b$ there will be exactly one positive root if $az - b > 0$, i.e. if $\frac{b}{a} < z$. Since our prime concern is with the large $z$ case, this condition does not significantly limit the applicability of this work. Thus, we have established that under the mild condition that $\frac{b}{a} < z$, $m_n$ is unimodal within the range $n > 0$, with the root of interest being the positive root in Eq. 3. As an aside we note also that for $|z| \gg |a|, |b|$, $\operatorname*{argmax}_{n} m_n \approx z$.

For convenience we operate on the log scale, denoting by $\Psi_n$ the log of the ratio of two consecutive terms in the series, i.e.

$$\Psi_n = \log(\omega_n) = \log\left(\frac{m_{n+1}}{m_n}\right) = \log\left(\frac{a + n}{b + n} \times \frac{z}{n + 1}\right) \qquad (4)$$

We can establish that the mode is a maximum by examining the derivative of $\Psi_n$ with respect to $n$.

$$\frac{d\Psi}{dn} = \frac{1}{a + n} - \frac{1}{b + n} - \frac{1}{n + 1}$$

Under the condition $a > 1$ or $a > b$ (and $b > 0$), this is negative for all $n$, and thus $\Psi_n$ is a strictly decreasing function of $n$. Since the log transformation is monotonic, $\omega_n$ must also be a strictly decreasing function of $n$. Therefore $\frac{m_{n+1}}{m_n}$ must be $> 1$ for $n < n_{\mathrm{m}}$ and $< 1$ for $n > n_{\mathrm{m}}$. Thus, function $m_n$ is indeed unimodal, and there is just one contiguous RoI. Thus, the CHF can be computed to arbitrary precision by considering a single contiguous range of terms centred on $n = n_{\mathrm{m}}$.

## 2.2. Finding the Limits of the Region of Interest

Our next step is to find a method for determining $n_{\mathrm{l}}, n_{\mathrm{u}}$, the lower and upper limits of the RoI for any required precision and pre-defined combination of the CHF parameters.

Formally, we seek $n_{\mathrm{l}}$ as the largest value of $n < n_m$ for which $\frac{m_n}{m_{n_{\mathrm{m}}}} \leq \varepsilon$, and $n_{\mathrm{u}}$ as the smallest value of $n > n_{\mathrm{m}}$ for which $\frac{m_n}{m_{n_{\mathrm{m}}}} \leq \varepsilon$. It is of note that this stopping criterion is considerably more conservative than the more conventional stopping criterion of requiring that the last term divided by the sum of all previous terms be below the tolerance [7].

### 2.2.1. Upper Limit

We will focus first on the upper end of the range, $n_\mathrm{u}$. Setting $n = n_\mathrm{m}$ in Eq. 2 and then applying it recursively yields the following expression:

$$\frac{m_{n_\mathrm{m}+k}}{m_{n_\mathrm{m}}} = \prod_{i=0}^{k-1} \frac{a+n_\mathrm{m}+i}{b+n_\mathrm{m}+i} \times \frac{z}{n_\mathrm{m}+i+1} \leq \varepsilon \qquad (5)$$

We denote $\Psi_n(k) = \log\left(\frac{m_{n+k}}{m_n}\right)$, so that $\Psi_n(1) \equiv \Psi_n$, with $\Psi_n$ defined in Eq. 4. The solution to Eq. 5 for $k$, $k_\mathrm{u}$, may be used to compute $n_\mathrm{u} = n_\mathrm{m} + k_\mathrm{u}$. We refer to $k_\mathrm{u}$ as the 'upper half-width' of the RoI. By working on the log scale, the log of a product becomes a sum of logarithms.

$$\Psi_{n_\mathrm{m}}(k) = \log\left(\prod_{i=0}^{k-1} \frac{a+n_\mathrm{m}+i}{b+n_\mathrm{m}+i} \times \frac{z}{n_\mathrm{m}+i+1}\right) = \sum_{i=0}^{k-1} \log\left(\frac{a+n_\mathrm{m}+i}{b+n_\mathrm{m}+i} \times \frac{z}{n_\mathrm{m}+i+1}\right)$$

$$\Psi_{n_\mathrm{m}}(k) = \sum_{i=0}^{k-1} \log(a+n_\mathrm{m}+i) - \sum_{i=0}^{k-1} \log(b+n_\mathrm{m}+i) + k\log(z) - \sum_{i=0}^{k-1} \log(n_\mathrm{m}+i+1)$$

By reverse application of the trapezoidal rule[21] with equi-spaced knots, we can replace the summation by an integration over $i$. The error bounds associated with the trapezoidal rule has a conservative upper limit of $O(m^{-2})$ where $m$ is the number of trapezoids[22], here equal to $k$, so it is to be expected that this approximation will be highly accurate even for moderate values of $k$.

$$\Psi_{n_\mathrm{m}}(k) \approx \int_{i=-1/2}^{k-1/2} \log(a+n_\mathrm{m}+i)\,di - \int_{i=-1/2}^{k-1/2} \log(b+n_\mathrm{m}+i)\,di + k\log(z) - \int_{i=-1/2}^{k-1/2} \log(n_\mathrm{m}+i+1)\,di$$

Having converted to an integral, we may conceptually regard $k$ as continuous rather than integer valued, as was done for $n$ in calculating the mode of $m_n$, in which case we may replace $\frac{m_{n_\mathrm{m}+k}}{m_{n_\mathrm{m}}} \leq \varepsilon$ by $\frac{m_{n_\mathrm{m}+k}}{m_{n_\mathrm{m}}} = \varepsilon$. Note that in keeping with the trapezoidal rule, the integral which best approximates the summation from $i = 0 \cdots k-1$, has limits from $-½$ to $k - ½$. However, for reasons of mathematical simplicity, we choose instead to integrate over the range $i = 0$ to $k$. For $n_m$ large, this will entail minimal loss of accuracy.

$$\log(\varepsilon) \approx \int_{i=0}^{k} \log(a+n_\mathrm{m}+i)\,di - \int_{i=0}^{k} \log(b+n_\mathrm{m}+i)\,di + k\log(z) - \int_{i=0}^{k} \log(n_\mathrm{m}+i+1)\,di \qquad (6)$$

The three logarithms in Eq. 6 can be integrated using the integration rule $\int \log(x+a)\,dx = (x+a)\log(x+a) - x$. However, there is no need to go through the process of performing this integration, because, as we will demonstrate, the result plays no role in the final expression for the limits of the RoI.

It is sufficient to note that using the expression above for $\int \log(x + a).dx$ we have derived an analytic expression for the logarithm of the target precision in terms of the number of summation terms required to reach the upper edge of the RoI, $k_u$. However, finding an analytic expression for the limit of the RoI necessitates inverting this and expressing $k_u$ in terms of the target precision. To this end, we use the Taylor's series expansion of $\log(\varepsilon)$, and assess increasing orders of the resulting polynomial in $k$ for accuracy in estimating the target function, $\log\left(\frac{m_{n_m+k}}{m_{n_m}}\right)$. We commence with a simple 1$^{st}$ order Taylor's polynomial:

$$\log(\varepsilon) \approx \Psi_{n_m}(0) + k \frac{d\Psi_{n_m}}{dk}(0)$$

Recalling that $\frac{a+n_m}{b+n_m} \times \frac{z}{n_m+1} = 1$, $\frac{d\Psi_{n_m}}{dk}(0) = 0$, so the coefficient of the linear term in the Taylor's expansion is zero. Further, $\Psi_{n_m}(0) = 0$, so we can turn immediately to a 2$^{nd}$ order Taylor's approximation.

$$\log(\varepsilon) \approx \frac{k^2}{2!} \frac{d^2 \Psi_{n_m}}{dk^2}(0)$$

From the expression for $\log(\varepsilon)$ in integral form above (Eq. 6), the expression for $\frac{d\Psi_{n_m}}{dk}$ is automatic, as it consists simply of differentiating an integral, which retrieves the original integrand, evaluated at $k$.

$$\frac{d\Psi_{n_m}}{dk}(k) = \log(a + n_m + k) - \log(b + n_m + k) + \log(z) - \log(n_m + 1 + k)$$

The second derivate is then given by

$$\frac{d^2 \Psi_{n_m}}{dk^2} = \frac{1}{a + n_m + k} - \frac{1}{b + n_m + k} - \frac{1}{n_m + 1 + k}$$

thus

$$\log(\varepsilon) \approx \frac{k^2}{2} \left( \frac{1}{a + n_m} - \frac{1}{b + n_m} - \frac{1}{n_m + 1} \right)$$

Our 2$^{nd}$ order Taylor's series approximation for the log precision thus yields the following approximation for the upper limit of the RoI.

$$k_u \approx \sqrt{\frac{2\log(\varepsilon)}{\left(\frac{1}{a + n_m} - \frac{1}{b + n_m} - \frac{1}{n_m + 1}\right)}} \quad (7)$$

It is of interest to consider higher order Taylor's polynomials as well. We follow the same procedure to extend the above method to include a 3$^{rd}$ order.

$$\log(\varepsilon) \approx \frac{k^2}{2}\left(\frac{1}{a+n_m} - \frac{1}{b+n_m} - \frac{1}{n_m+1}\right) - \frac{k^3}{6}\left(\frac{1}{(a+n_m)^2} - \frac{1}{(b+n_m)^2} - \frac{1}{(n_m+1)^2}\right) = g(k) \qquad (8)$$

The solution to Eq. 8 remains tractable, requiring nothing more than a standard method for solving cubics [23]. In the event that the cubic $g(k) - \log(\varepsilon)$ has 3 real roots, we must have a method of determining which is the root of interest. For the form of cubic presented in Eq. 8, this is a straightforward exercise. The absence of a linear term means that the cubic must have two turning points. Denoting the cubic and quadratic coefficients by $C_3$ and $C_2$ respectively, $g(k) = C_2 k^2 + C_3 k^3$, these turning points are at $k = 0$ (maximum, $\log(\varepsilon) = 0$) and $k = -\frac{2C_2}{3C_3}$ (minimum). From inspection of the definitions of $C_3$ and $C_2$, we can see that under the condition $a \geq 1$, $C_3 > 0$ and $C_2 < 0$, and thus $-\frac{2C_2}{3C_3} > 0$. Note that the turning point at $-\frac{2C_2}{3C_3}$ exists only in the Taylor's approximation to $\log\left(\frac{m_{n_m+k}}{m_{n_m}}\right)$. The unapproximated function $\log(\varepsilon)$, as we have seen, is unimodal. If there are 3 roots, then they must be such that one is for $k < 0$, one is for $0 < k < -\frac{2C_2}{3C_3}$ and the third is for $k > \frac{2C_2}{3C_3}$. The third of these is not of interest because at this point the Taylor's polynomial $g(k)$ has already lost fidelity with the function $\log\left(\frac{m_{n_m+k}}{m_{n_m}}\right)$ which it seeks to approximate. Therefore, the root of interest is that which lies in the range $0 < k < -\frac{2C_2}{3C_3}$.

Since $C_3 > 0$, $\lim_{k \to \infty} g(k) = \infty$ and $\lim_{k \to -\infty} g(k) = -\infty$. This means that if there is only one root to the cubic, it must be for $k < 0$, and this is not of interest.

The point where $k = -\frac{2C_2}{3C_3}$, being a local minimum of $g(k)$, also represents the minimum value of the precision that can be estimated for any combination of parameter values. Beyond this point $g(k)$ starts to increase, whereas the function $\Psi_{n_m}(k) = \log\left(\frac{m_{n_m+k}}{m_{n_m}}\right)$ which it seeks to approximate continues to decrease. It is of interest to express this minimum precision value ($\log(\varepsilon_{\min})$, the minimum value of $\log(\varepsilon)$ for any value of $k > 0$, which is the most precise estimate for which the RoI limits can be estimated) as a function of the CHF parameters. Fortunately, this is simply a matter of setting $k = -\frac{2C_2}{3C_3}$ in $g(k)$. In terms of the cubic and quadratic coefficients $C_3$ and $C_2$ defined as above, $\log(\varepsilon_{\min}) = \frac{4C_2^3}{27C_3^2}$. We will see in the Results section that for a wide range of CHF parameter values, this yields a precision well in excess of that required for most practical applications (smaller $\log(\varepsilon_{min})$ than necessary).

We thus have two estimates for the number of summation terms required to estimate the CHF to pre-specified precision $\varepsilon$, the first given by Eq. 7 based on a 2nd order Taylor approximation to Eq. 6, and the 2$^{nd}$ by the solution to the 3$^{rd}$ order Taylor approximation Eq. 8. The accuracy of these approximations will be assessed and compared in the Results section.

### 2.2.2. Lower limit

Estimation of the lower limit of the RoI, $n_l$ proceeds analogously to estimation of the upper limit. We commence by inverting the forward recurrence relation in Eq. 1, to yield a backward recurrence relation, $\frac{m_{n-1}}{m_n} = \frac{b+n-1}{a+n-1} \times \frac{n}{z}$. From that starting point, the same steps can be followed as for the upper limit presented above. Due to its similarity to the derivation of the upper limit, we do not present the derivation in full. Instead, we observe that if $\Psi_n^{-1}(k) = \log\left(\frac{m_{n-k}}{m_n}\right)$ is the recurrence relation in the backwards direction, then $\Psi_n^{-1}(k) = \Psi_n(-k)$. This means that the 2$^{nd}$ order Taylor's approximation for the lower half-width, which is dependent on $k$ only through $k^2$, is identical to that for the upper half-width. Further, the expression for the precision in the lower direction when using a 3$^{rd}$ order Taylor's approximation differs from that of the upper direction only by an inversion of the sign for the cubic term. It is presented in Eq. 9 for completeness.

$$\log(\varepsilon) = \frac{k^2}{2}\left(\frac{1}{a+n_m} - \frac{1}{b+n_m} - \frac{1}{n_m - 1}\right) + \frac{k^3}{6}\left(\frac{1}{(a+n_m)^2} - \frac{1}{(b+n_m)^2} - \frac{1}{(n_m - 1)^2}\right) \quad (9)$$

The limits of the parameters at which the Taylor's approximation breaks down will be explored in the Results section.

For calculation of the upper half-width, the signs of the coefficients of the Taylor's series approximations alternate ($C_2 < 0$, $C_3 > 0$, $C_4$ will be $< 0$, etc.), implying a sequence with terms of alternating sign and hence oscillating convergence. Oscillating convergence is often very slow, but there exist a variety of methods of speeding it up. Many modern acceleration methods are summarized by Bender and Orszag [24] and Cohen et al [25]. A complicating consideration in this case is that the elements of the sequence to be summed are consecutive terms in a Taylor's polynomial, the highest power of which must be kept low to ensure that the half-width, $k_u$, which is a root of the polynomial, can be readily determined computationally. This is required because a root finding method needs to be applied to the resulting polynomial to achieve our goal of estimating the number of terms in the RoI of the CHF. We thus restrict ourselves to perhaps the best known of these acceleration methods, Euler's transformation, and straightforward variations thereof which can be applied even for very short series.

One method to generate Euler's transformation of an alternating sequence $(-1)^k b_k$ is to consider the sequence of partial sums $a_n = \sum_{k=0}^{n}(-1)^k b_k$ and iteratively replace each pair of consecutive terms in $a_n$ by its mean, $a'_n = \frac{1}{2}(a_n + a_{n+1})$, and then take the last terms of the resulting derived sequences, $a'_1, a''_2, \cdots$.[26] Euler's summation method is known to significantly accelerate the convergence of alternating (oscillating) series. Van Wijngaarden noted that very often convergence is optimised by not following this process through to the end, but terminating it $\frac{2}{3}$ of the way through [27]. Arguing intuitively, when the damping rate of the consecutive terms is very gradual ($|b_{k+1}|$ only slightly smaller than $|b_k|$), performing the averaging process once only, which is equivalent to taking the original series with just half of the final term, will achieve good acceleration of the convergence. We thus satisfy ourselves with just such a truncated version of Euler's summation. This is tantamount to summing the original series, but including only half of the final term. We thus have in total four Taylor's series approximations with varying numbers of terms for consideration – two with a whole number of terms (quadratic and cubic), and another one corresponding to each of these with half of the last term deleted. We will denote these four series as T1.5, T2, T2.5 and T3 in the Results section. We note that it is only in the upper direction that the convergence is alternating, and thus a convergence acceleration method is needed in the upper direction only. However, to be consistent, we experiment with all four methods in both the upper and lower directions.

In summary, the proposed method of estimating the CHF consists of the following steps:

1. Calculate the centre of the RoI, $n_m$, using Eq. 3
2. Estimate the upper and lower limits of the RoI from the roots of one of the Taylor's polynomials T1.5, T2 (Eq. 7), T2.5 or T3 (Eqs. 8 and 9)
3. Calculate the value of $m_n$ at the lower edge of the RoI, $m_{n_l} = \frac{a^{(n_l)} z^{n_l}}{b^{(n_l)} n_l!}$
4. Use the recurrence relation in Eq. 2 to calculate all subsequent values of $m_n$ until the upper edge of the RoI, $m_{n_u}$.
5. Sum all the component terms, $M(a, b, z) = \sum_{n=n_l}^{n_u} m_n$

Steps 1 and 2 (in particular step 2) constitute the novelty of this work and have been described in detail above. The easiest way to implement step 3 is by noting the relationship between the rising factorial and the Gamma function as $x^{(n)} = \frac{\Gamma(x+n)}{\Gamma(x)}$. However, for large values of $a$, $b$, or $n_l$, it may be necessary to be more judicious with the order of computation to avoid overflow or underflow in intermediate calculations. Steps 4 and 5 are then straightforward.

## 3. Results

Given that the novelty of this work lies in steps 1 and 2 above, the aim of the first set of results is to corroborate the proposed method for these steps by way of simulation.

Step 1, calculating the centre of the RoI as the solution to a quadratic equation (Eq. 3), being purely analytic, requires no experimental verification. Step 2, however, involves approximations in the form of replacing a summation by an integral (inverse trapezoidal rule) and then in the form of replacing the resulting function by its Taylor's series expansion. Our first set of experiments thus explores the accuracy of these approximations for a wide variety of CHF parameter settings.

For purposes of comparison, the exact relationship between the number of terms included in the calculation of the CHF on either side of the maximum point, $n_m$, and the resulting precision achieved can be determined as follows. Using an exhaustive 'increment-and-check' approach, the recurrence relation in Eq. 2 is applied iteratively starting at $n = n_m$ in the upward direction, checking for convergence to within pre-determined precision after each term is added. This can then be repeated using the reverse recurrence relation for the downward direction. (This avoids inaccuracies due to adding numbers which differ by many orders of magnitude, which would arise were the summation to commence at $n = 0$. The limitations of commencing summation at $n = 0$, as is usually done when applying series-based methods for calculating the CHF, is explored more fully in section 3.1 below.) This can then be compared to the precision calculated using either Eq. 6a for the quadratic approximation to the RoI width, or Eqs. 8 and 9 for the cubic approximation.

Results are shown in Figures 1 and 2 which show the index number of the lowest and highest term included in the RoI (x-axis) and the corresponding precision (y-axis). The precisions shown are the actual precision ('Exact'), which uses the increment-and-check approach, and the precision estimated using the inverse trapezoidal rule plus Taylor's expansion, considering the various polynomial powers discussed in the Methods section. Since the CHF with $z < 0$ can be transformed to a CHF with $z > 0$ using Kummer's transformation, it is sufficient to consider the case $z > 0$ only. We focus on the case where $z$ is large, since the small $z$ case is easily solved by directly applying the recurrence relation in Eq. 2 starting at $n = 0$ without pre-calculating an RoI. It has been shown[11] that the recurrence relation works very well until at least $z = 50$, so in our experiments we consider only values of $z$ larger than this. Figure 1 shows results for scale parameter $z = 50$, and Figure 2 shows results for $z = 100$. A wide range of values of $a$ and $b$ are considered across the panel rows and columns respectively.

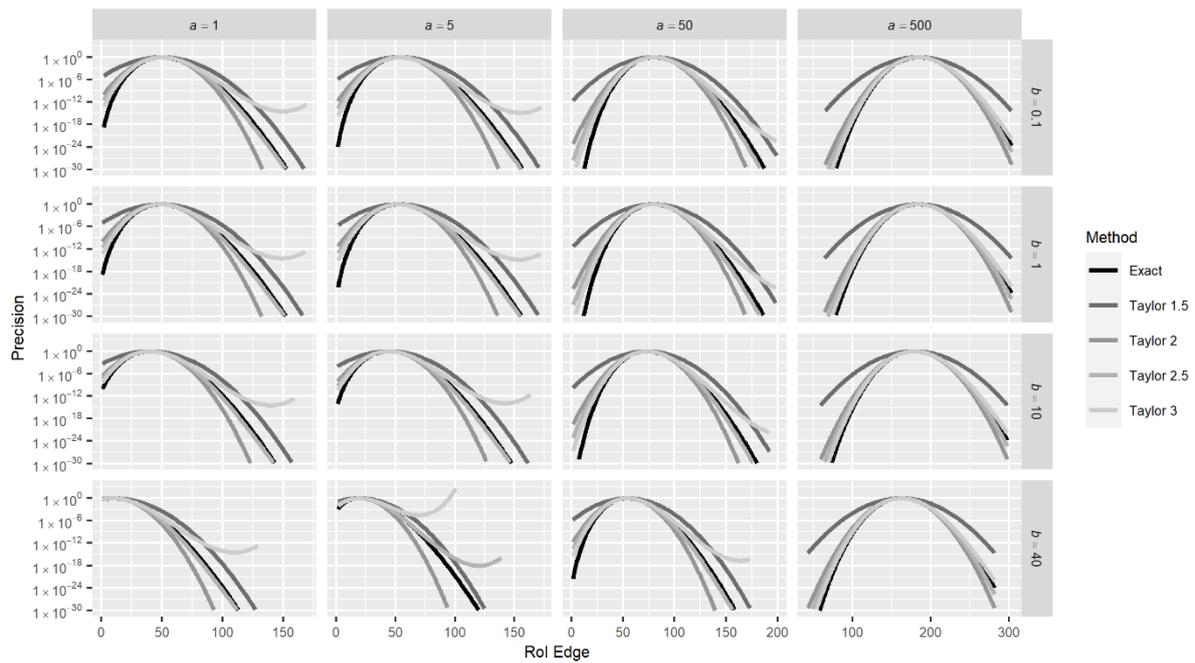

**Fig. 1** Approximations of the precision achieved using each of four Taylor's polynomials as a function of the RoI edge, with scale parameter $z = 50$, compared with the exact precision achieved

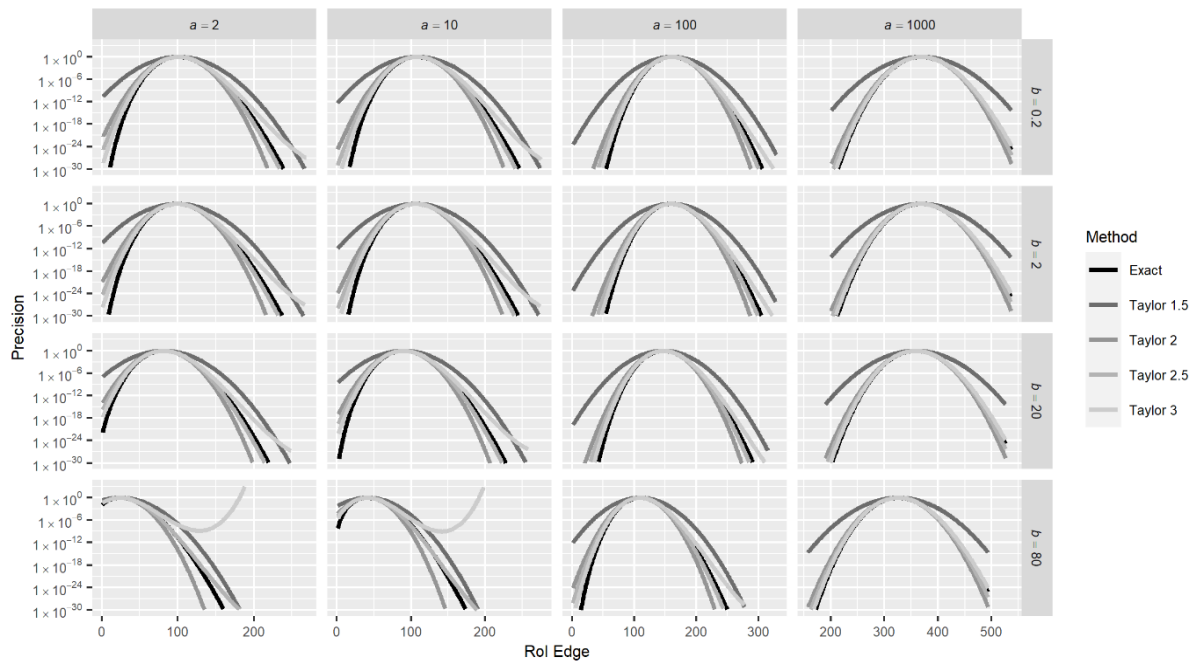

**Fig. 2** Approximations of the precision achieved using each of four Taylor's polynomials as a function of the RoI edge, with scale parameter $z = 100$, compared with the exact precision

The x-axis in Figures 1 and 2 is the RoI edge. Values above the mode represent the upper edge of the RoI and values below the mode represent the lower edge of the RoI. The y-axis represents the ratio of the CHF component term at the RoI edge to the maximum term. (Thus, the mode of each curve is equal to 1 by definition.) Columns represent different values of $a$ ranging from much smaller than $z$ to much larger than $z$.

In results not shown it was found that the Taylor's approximation breaks down for $b$ greater than both $z$ and $a$, and thus the RoI method is not appropriate for this condition. However, this is not a significant limitation because $b$ only appears in the denominator of $M(a, b, z)$ and thus, for large $b$, standard non-RoI approaches for calculating $M(a, b, z)$ will converge quickly. The formula provided in Eq. 3 for $n_\text{m}$ corroborates this, and similar findings have been reported in the literature[7]. Thus, rows in Figures 1 and 2 represent a range of values of $b$ from very small to slightly below $z$.

Figures 1 and 2 show results for $a > 0$ and $b > 0$ only. However, the approximations work equally well for $b < 0$. Results are shown in the Appendix. For $a < 0$, the quadratic formula in Eq. 3 often has a negative discriminant and thus no real solution, so methods other than the RoI need to be used.

The first overarching point to note from Figures 1 and 2 is just how accurately the 'RoI vs. precision' relationship can be characterised using a Taylor's approximation with as few as two terms. For either $a$ or $z$ substantially larger than $b$, T2, T2.5 and T3 all yield precision estimates very close to that of the exact precision over a wide range of RoI limits. T1.5 is not as accurate but is always conservative, meaning it overestimates the precision (i.e. underestimates $\varepsilon$) for any given RoI width. T2.5 is generally more accurate than T3.

Both T2.5 and T3 must eventually break down as the required precision becomes very small since they are cubic polynomials with a positive coefficient on the cubic term and hence must approach infinity as the RoI upper edge becomes sufficiently large. The exact precision, on the other hand, being unimodal, always approaches zero for large RoI upper edge. This divergence between the exact and estimated precisions can be observed for T2.5 and T3 in Figure 1 and for T3 in Figure 2, when $a$ is small and $b$ approaches $z$.

It is thus of interest to explore the range of CHF parameters for which this method can be utilised for varying values of precision. To this end we use our previously established relationship $\log(\varepsilon_\text{min}) = \frac{4C_2^3}{27C_3^2}$ to characterise $\varepsilon_\text{min}$, the minimum value of $\varepsilon$ for which RoI limits are calculable for a given combination of the CHF parameters. $C_2$ and $C_3$ are the quadratic and cubic coefficients in T2.5. Results are presented in Figures 3 and 4.

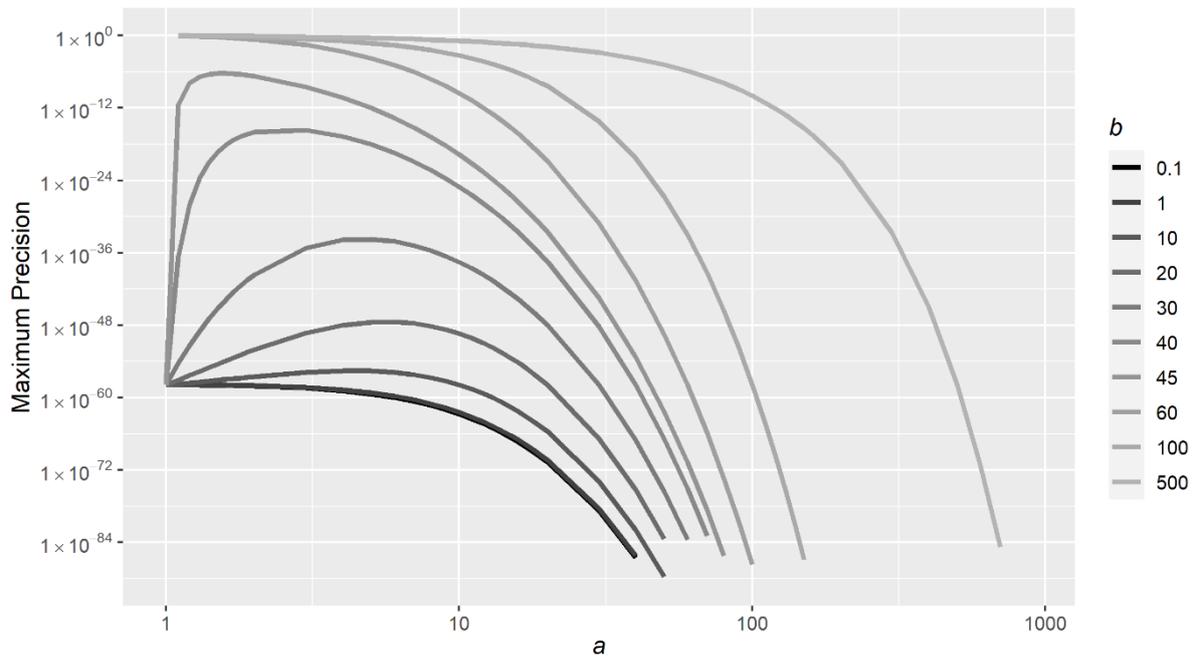

**Fig. 3** Maximum achievable precision ($\varepsilon_{min}$) as a function of $a$ and $b$, for $z = 50$ using T2.5

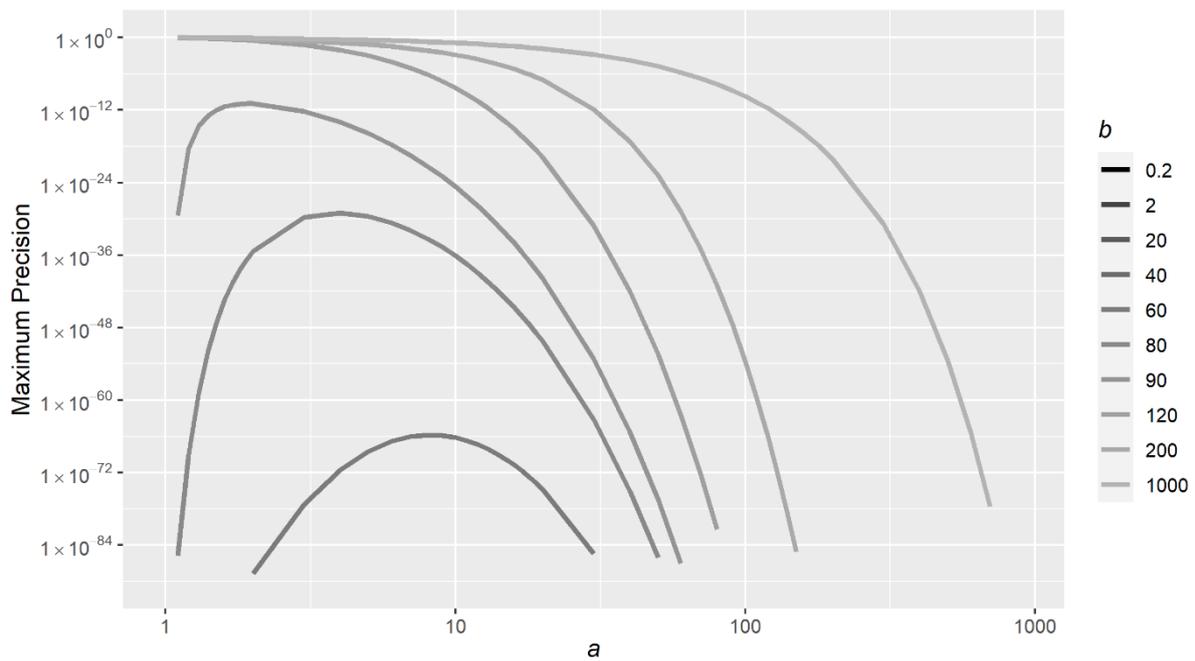

**Fig. 4** Maximum achievable precision ($\varepsilon_{min}$) as a function of $a$ and $b$, for $z = 100$ using T2.5

Figures 3 and 4 show that the Taylor's approximation with 2.5 terms can be used to estimate the required RoI width for pre-determined required precision up to very small precisions, except for situations where both $b \gtrsim z$ and $a \ll z$. With $z = 50$, T2.5 can be used to estimate the required RoI limits for a required precision of $10^{-6}$ or smaller for any value of $a$. Note that this is not the precision of the RoI width itself, but the maximum precision (smallest $\varepsilon$) for which the T2.5 approximation will be able to produce an estimate for the RoI upper edge. Figures 1 and 2 show that as $\varepsilon$ approaches this smallest value, the estimate of the RoI edge becomes more imprecise. However, they also show that it is always conservative, guaranteeing that the achieved precision will be at least as good as that planned. As $z$ increases, the range of applicability of this method expands to ever smaller precisions.

## 3.1. Efficiency

The computational efficiency advantage of the RoI method over traditional series summation methods for calculating the CHF is manifest in two ways. Firstly, the number of component terms that need to be included in the summation is reduced. Secondly, a priori calculation of the boundaries of the RoI obviates the need for an 'increment-and-check' algorithm to check for convergence after addition of each term or at least periodically, which would necessitate additional flow control logic and also impedes the use of vectorisation in the algorithmic implementation. It is likely that the efficiency advantage of the second benefit is platform dependent, but the advantage of the first is easily quantified by comparison of the number of summation terms required by each method. The number of summation terms required when using the RoI method is simply given by $n_u - n_l$, with both $n_u$ and $n_l$ calculated using the chosen root-finding approximation method. The increment-and-check method starts at term $m_0$, and continues until the calculated term $m_n$ is less than the maximum by a factor equal to the required precision. Figure 5 compares the number of summation terms required as estimated by the RoI method using T2.5 to the number required using the conventional 'increment-and-check' approach (IC) starting from $n = 0$ without pre-calculating an RoI, to achieve precision levels of $10^{-6}$, $10^{-12}$ and $10^{-18}$.

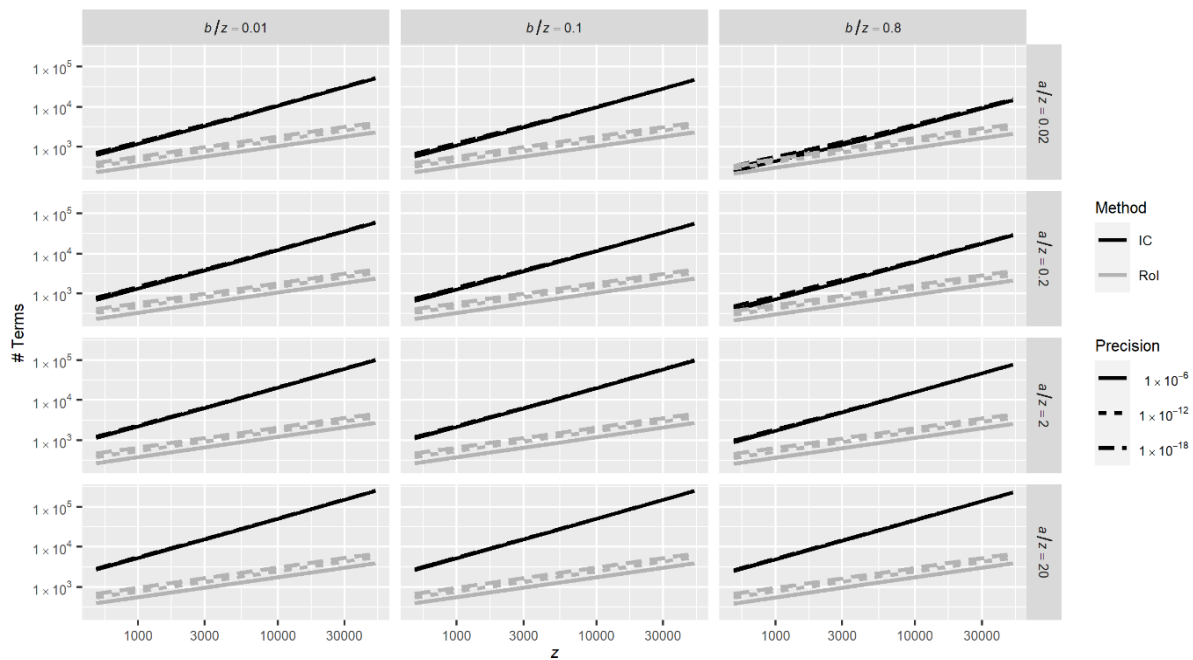

**Fig. 5** Number of summation terms required to achieve desired precision using both the RoI method (T2.5) and the 'increment-and-check' method which starts at $n = 0$ (IC)

Figure 5 shows that the number of summation terms required using the RoI approach is reduced by a factor of 10 to 100 relative to the number required using the IC approach over the range of parameter values considered. The reduction factor is dependent primarily on the value of $z$, with a reduction factor of ~10 being achieved for $z$ in the low thousands, and a reduction factor of ~100 being realised for $z$ in the vicinity of 10,000.

## 3.2. Precision

Aside from inefficiency, another problem with increment-and-check approaches to calculation of the CHF are loss of accuracy due to multiplication of numbers across vastly different orders of magnitude which ultimately lead to computational overflow (or underflow). The exact point at which the increment-and-check approaches breaks down will depend primarily on the number of bits of precision available for floating point calculation and also on the order of operations. In order to present a platform independent representation of the breakdown point of the increment-and-check approach, we use the range of intermediate terms (ratio of largest to smallest) as a surrogate for the loss of accuracy and likelihood of overflow. The larger this ratio, the more likely it is that the increment-and-check approach will suffer from insufficient accuracy, and ultimately, fail.

Our interest is primarily in the case where the scale parameter $z$, is large, and under this condition the value of the CHF rapidly increases to the point of computational overflow. To circumvent this, we choose to work with a related function, a discrete valued statistical distribution known as the Poisson-Beta (PB) distribution, which, because it is a statistical distribution, is guaranteed to return values <1 for any combination of its parameters. Representing the variate of the PB distribution by $x$, the relationship between the PB distribution $f_{PB}$ and the CHF ($M$) is given by the following:

$$f_{PB}(x;\alpha,\beta,\gamma) = \frac{\gamma^x}{x!} \frac{(\alpha)^{(x)}}{(\alpha+\beta)^{(x)}} e^{-\gamma} M(\beta, \alpha+\beta+x, \gamma)$$

We adhere to the conventional notation for the parameters of the PB distribution, $\alpha, \beta, \gamma$, and thus note that the mapping of the parameters of the PB distribution to the arguments of the CHF is $a \leftarrow \beta, b \leftarrow \alpha + \beta + x$ and $z \leftarrow \gamma$. In producing Figure 6, for simplicity we set $x = 0$. The required precision, $\varepsilon$, was set to $10^{-12}$.

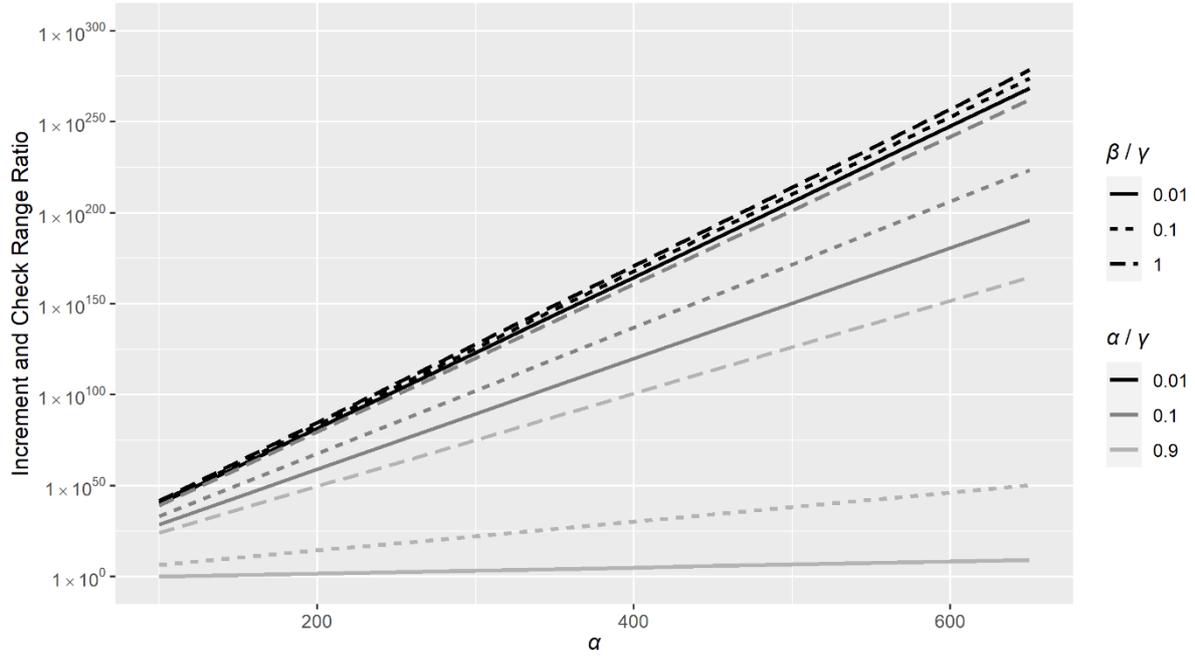

**Fig. 6** Ratio of the smallest to the largest term in the summation using the increment-and-check approach to calculating the density of the PB distribution

It was not possible to calculate the PB distribution density for values of $\gamma$ above ~700 (x-axis) because beyond this point the extremely high ratios between the largest and smallest terms in the summation caused either computational overflow or computational underflow, depending on the order of the operations.

In contrast, no such limit was found with the RoI method. The ratio of the largest to the smallest term in the summation was <$10^{21}$ for all parameter settings considered, and was generally around $10^{12}$, as expected for a planned precision of $10^{-12}$.

In contrast, using the RoI method, no upper limit for $\gamma$ could be found beyond which the PB distribution density could no longer be calculated accurately. Accuracy tests for the RoI method were undertaken by calculating $\sum_x f_{PB}(x; \alpha, \beta, \gamma)$, which, since $f_{PB}$ is a density function, must sum to 1 when the range of summation $x$ is wide enough. Tests were performed with values of $\gamma$ as high as 200,000, and accuracy was still within $10^{-9}$.

## 4. Discussion

The 'RoI' method presented in this work massively extends the range of applicability and enhances the efficiency of the recurrence relation based methods for computing the CHF, the simplest methods available for computing this function. Such recurrence relation-based methods, have, until now, been limited to small to moderate values of the CHF scale parameter. Using the RoI method, no upper limit for the scale parameter beyond which the method could not be applied, could be found. We have proven that the component terms of interest needed to calculate the CHF accurately constitute a single contiguous range and have provided a method for accurately estimating the upper and lower limits of that range by nothing more than solving for a cubic whose coefficients are easily calculable functions of the CHF parameters.

The estimation of the RoI limits becomes increasingly precise as $z$ increases. Results shown are for values of $z$ low enough for simpler 'increment-and-check' recurrence relation based approaches could still feasibly be used, and thus show the RoI approach at its weakest. Even under these conditions, the RoI approach has been shown to be robust.

## 5. Conclusion

When computing the Confluent Hypergeometric Function $_1F_1(a;b;z) = \sum_{n=0}^{\infty} \frac{a^{(n)}z^n}{b^{(n)}n!}$, with $a, b, z \in \mathbb{R}$, a 'region of interest' on $n$, outside of which the component terms $\frac{a^{(n)}z^n}{b^{(n)}n!}$ make negligible contribution to the summation, can be easily identified as the solution to a quadratic or cubic equation whose coefficients are directly calculable functions of $a, b, z$. This massive extends the applicability of recurrence relation based approaches to computing $_1F_1(a;b;z)$, until now restricted to small values of $z$ ($\lesssim 100$), to include arbitrarily large $z$.

## 6. Appendix

Figures 7 and 8 show the T1.5, T2, T2.5 and T3 approximations to the exact precision achieved for ranges of values of $b < 0$ and $a > 0$, for $z = 100$ (Figure 7) and $z = 100$ (Figure 8). Only non-integer values of $b$ are considered in order to avoid the division by zero problem that occurs for $b \in \mathbb{Z}^-$. As for positive values of $b$, the T2.5 approximation to the exact precision achieved shows exceptional fidelity for a wide range of values of $a$ and $b$, even for moderate values of $z$.

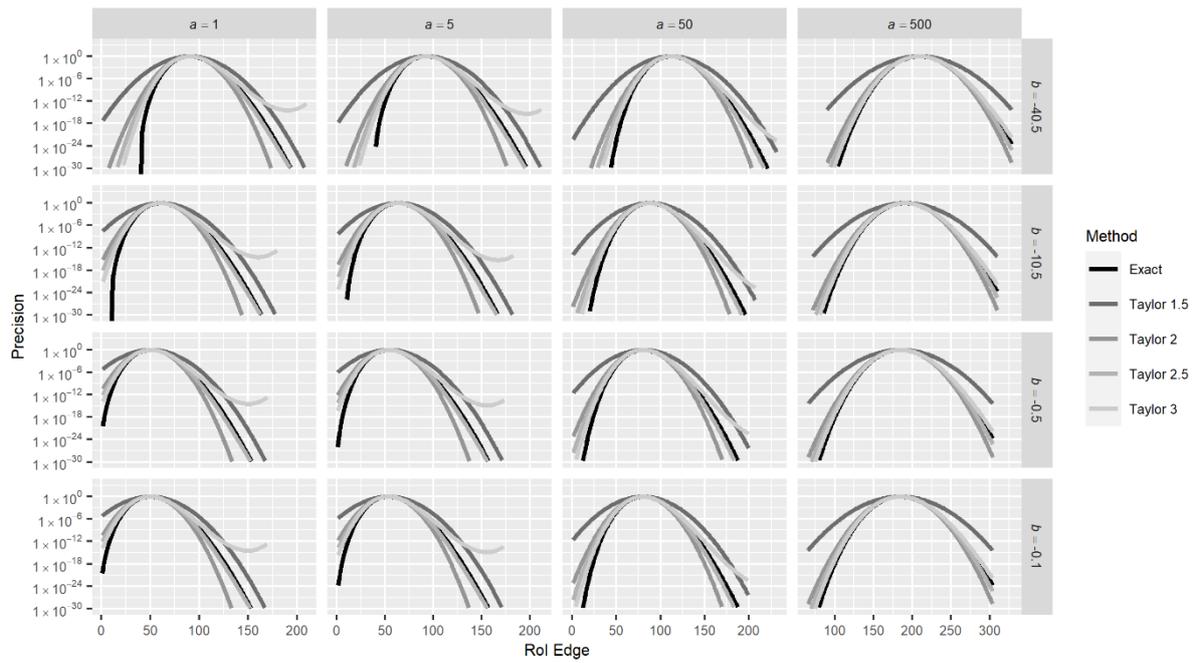

**Fig. 7** Approximations of the precision achieved using each of four Taylor's polynomials as a function of the RoI edge, with scale parameter $z = 50$, considering negative values of $b$, compared with the exact precision achieved

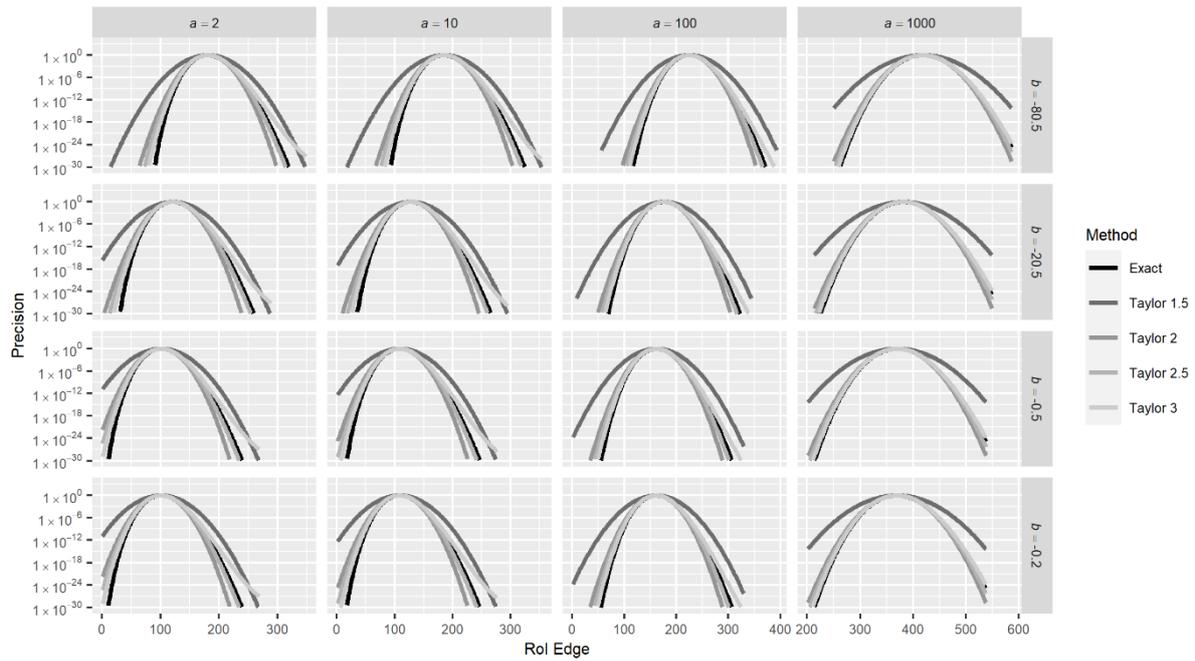

**Fig. 8** Approximations of the precision achieved using each of four Taylor's polynomials as a function of the RoI edge, with scale parameter $z = 100$, considering negative values of $b$, compared with the exact precision achieved

# 7. Declarations

**Ethical Approval:** This study did not involve any human or animal subjects, nor did it involve disclosure of any private data. Thus, consent to participate and consent to publish are not applicable.

**Availability of supporting data:** Data sharing is not applicable to this article as no datasets were generated or analysed during the current study.

**Competing interests:** The author has no competing interests to declare that are relevant to the content of this article.

**Funding:** A contribution towards this work in the form of a seed funding grant, was made by AusTriM, the Australian Trials Methodology group.

**Authors' contributions:** As sole author, Alan Herschtal was responsible for development of the methodology, coding, write-up of results and review of the paper for submission. Other contributors are duly acknowledged below.

**Acknowledements:** The encouragement and guidance of Stephane Heritier, Monash University, School of Public Health and Preventive Medicine, is acknowledged. Proof reading by Liz Ryan, Monash University, School of Public Health and Preventive Medicine, is acknowledged. The support of AusTriM, the Australian Trials Methodology group, is acknowledged.